\documentclass{article}
\usepackage{latexsym}
\usepackage{amssymb}
\usepackage{amsthm}
\usepackage{amsmath,amsfonts,mathabx}
\usepackage[all]{xy}

\newtheorem{theorem}{Theorem}
\newtheorem{proposition}{Proposition}

\newtheorem{definition}{Definition}
\newtheorem{corollary}[theorem]{Corollary}

\begin{document}
\title{\textbf{Orbifolds, geometric structures and  foliations. Applications to harmonic maps}}
\author {Robert A. Wolak \\{\small Wydzial Mathematyki i Informatyki, Uniwersytet Jagiellonski}\\{\small Lojasiewicza 6,  30-348 Krakow, Poland, robert.wolak@im.uj.edu.pl.}\\
\\ Dedicated to Anna Maria Pastore on her 70th birthday
}

\maketitle

\vskip 2mm
{\footnotesize\noindent{\em$2011$ Mathematics Subject Classification\/}.  53C12 }


\begin{abstract}
In recent years a lot of attention has been paid to topological spaces which are a bit more general than smooth manifolds - orbifolds. Orbifolds are intuitively speaking manifolds with some singularities. The formal definition is also modelled on that of manifolds, an orbifold is a topological space which locally 
is homeomorphic to the orbit space of a finite group acting on $R^n$. Orbifolds were defined by Satake, as V-manifolds, cf. \cite{Sa},  then studied by W. Thurston,  cf. \cite{Th}, who introduced the term "orbifold". 
Due to their importance in physics, and in particular in the string theory, orbifolds have been drawing  more  and more attention. In this paper we propose to show that the classical theory of geometrical structures, cf. \cite{St,Fu,Mi}, easily translates itself to 
the context of orbifolds and is closely related to the theory of foliated geometrical structures, cf. \cite{Wo0}. Finally, we propose a foliated approach to the study of harmonic maps between Riemannian orbifolds based on our previous research into transversely harmonic maps, cf. \cite{KW1,KW2}.

\end{abstract}
\vskip 2mm

\section{Orbifolds and their smooth complete maps}


 Let $X$ be a topological space, and fix $n \geq 0.$  An $n$-dimensional orbifold chart on $X$ is given by a connected open subset $ \tilde{U}\subset  R^n,$  a finite group $\Gamma$ of smooth diffeomorphisms of $\tilde{U}$, 
and a map $\phi \colon \tilde{U} \rightarrow X$ so that $\phi$  is $\Gamma$-invariant and induces a homeomorphism of  $\tilde{U}/{\Gamma} $ onto an open subset $U \subset X.$

\medskip

An embedding $ \lambda \colon  (\tilde{U}, \Gamma ,\phi ) \rightarrow  (\tilde{V} , \Delta ,\psi )$ between two such charts is a
smooth embedding $ \lambda  \colon \tilde{U}\rightarrow  \tilde{V}$  with $ \psi \lambda  = \phi .$ Then there exists a homomorphism $\alpha \colon \Gamma \rightarrow \Delta$ such that the mapping $\lambda$ is $\alpha$-equivariant. 

\medskip

 An orbifold atlas on $ X$ is a family ${\mathcal U } = \{(\tilde{U}, \Gamma ,\phi )  \}$ of such charts, which cover $X$ and are locally compatible:

\medskip

\noindent
 given any two charts  $(\tilde{U}, \Gamma ,\phi ) $  for $U = \phi (\tilde{U}) \subset X$  and  $(\tilde{V}, \Delta ,\psi ) $  for $V = \psi (\tilde{V}) \subset X$
 and a point $x \in  U \cap V , $ there exists an open neighborhood $ W \subset  U \cap V $ of $x$ and a chart   $(\tilde{W}, \Lambda ,\tau ) $  for $W = \tau (\tilde{W}) \subset X$
such that there are two embeddings  $(\tilde{W}, \Lambda ,\tau ) \rightarrow   (\tilde{U}, \Gamma ,\phi )$   and $(\tilde{W}, \Lambda ,\tau ) \rightarrow (\tilde{V}, \Delta ,\psi ).$ 

\medskip

An atlas ${\mathcal U }$  is said to refine another atlas ${\mathcal V }$  if for every chart in ${\mathcal U }$ there
exists an embedding into some chart of ${\mathcal V }.$ Two orbifold atlases are said to be equivalent if they have a common refinement.

\medskip

With these notions defined above we can formulate the definition of an orbifold

\begin{definition}  An effective orbifold $X$ of dimension n is a paracompact Hausdorff space $X$ equipped with an equivalence class $ [{\mathcal U}] $ of n-dimensional orbifold
atlases

\end{definition}

In this paper we will consider only effective orbifolds without further mentioning it. 

\medskip


For any effective orbifold we can find an atlas as in the following definition due to Borzellino and Brunsden, cf. \cite{B}

\begin{definition}
 An n-dimensional smooth orbifold $\mathcal O$ consists of a paracompact,
Hausdorff topological space $X_O$ called the underlying space, with the following
local structure. 
	For each $x \in X_O$ and neighborhood $U$ of $x,$ there is a neighborhood
$U_x \subset U$ of $x,$ an open set $\tilde{U}_x$ diffeomorphic to $R^n,$ a finite group $\Gamma_x$ acting smoothly and effectively on $\tilde{U}_x$ which fixes $0 \in \tilde{U}_x$, and a homeomorphism $\phi_x : \tilde{U}_x/\Gamma_x \rightarrow U_x$
with $\phi_x(0) = x.$ 
	These actions are subject to the condition that for a neighborhood
$U_z \subset U_x$ with the corresponding $\tilde{U}_z \tilde{=} R^n,$  group $\Gamma_z$ and homeomorphism $\phi_z \colon  \tilde{U}_z/\Gamma_z \rightarrow U_z$ there is a smooth embedding  $\tilde{\psi}_{zx} \colon \tilde{U}_z \rightarrow \tilde{U}_x$ and an injective homomorphism  $\theta_{zx}\colon \Gamma_z \rightarrow  \Gamma_x$ so that $\tilde{\psi}_{zx}$ is equivariant with respect to $\theta_{zx}.$ This means that for any $\gamma \in
\Gamma_z ; \tilde{\psi}_{zx}( \gamma\tilde{y}) = \theta_{zx}(\gamma ) \tilde{\psi}_{zx}(\tilde{y})$ for all $\tilde{y}\in \tilde{U}_z$ and  the following diagram commutes:

\begin{equation*}
\xymatrixcolsep{7pc}\xymatrix{ 
\tilde{U}_z \ar[r]^{\tilde{\psi}_{zx}} \ar[d]  &  \tilde{U}_x \ar[d] \\
\tilde{U}_z/\Gamma_z  \ar[r]^{{\psi}_{zx} = \tilde{\psi}_{zx}/\Gamma_z } \ar[d]^{\phi_z}    &  \tilde{U}_x/\theta_{zx}(\Gamma_z)  \ar[d]^{\phi_x}\\
{U}_z \ar[r]^{\subset}  & {U}_{x}
 \\}
\end{equation*}

\end{definition}

For a given orbifold atlas $\mathcal U$, the underlying topological space $X_O$ is homeomorphic to the topological space $X_{\mathcal U}$  defined as the quotient of the disjoint union of $\tilde{U}_i$ by the equivalence relation induced by the action of 
the groups $\Gamma_i$ and the compatibility condition. When there is no ambiguity the underlying space of the orbifold $X$ will be denoted by the same letter $X.$

We will refer to the neighborhood $U_x $ or $(\tilde{U}_x, \Gamma_x)$ or $ (\tilde{U}_x, \Gamma_x, \rho_x, \phi_x)$ as an orbifold chart. In the 4-tuple notation, we are making explicit the
representation $\rho_x : \Gamma_x \rightarrow  Diff(\tilde{U}_x).$ The isotropy group of $x$ is the group $\Gamma_x.$ The definition of orbifold implies that the germ of the action of $\Gamma_x$ in a neighborhood of the origin of $R^n$ is unique, so that by shrinking $\tilde{U}_x$ if necessary, $\Gamma_x$ is well-defined up to isomorphism. The singular set of O is the set of points $x \in X_O$ with $\Gamma_x \not= \{ e\}.$ 

\medskip

We present the definition of a complete orbifold map from \cite{B}.

\begin{definition}
 A $C^{\infty}$ complete orbifold map $(f, \{\tilde{ f_x}\}, \{\Theta_{f,x}\})$ between smooth orbifolds $\mathcal O$ and $\mathcal P$ consists of the following

(1) a continuous map $f \colon X_{\mathcal O} \rightarrow X_{\mathcal P}$ of the underlying topological spaces,

(2) for each $y \in X_{\mathcal O},$ there is  a group homomorphism $\Theta_{f,y} \colon \Gamma_y \rightarrow \Gamma_{f(y)},$

(3) a smooth $\Theta_{f,y}$-equivariant lift $\tilde{f}_y\colon \tilde{U}_y \rightarrow
\tilde{V}_{f(y)}$  where $(\tilde{U}_y, \Gamma_y)$ is an orbifold
chart at $y$ and $( \tilde{V}_{f(y)}, \Gamma_{f(y)})$ is an orbifold chart at $f(y). $
That is the following diagram commutes:

\begin{equation*}
\xymatrixcolsep{7pc}\xymatrix{ 
\tilde{U}_y \ar[r]^{\tilde{f}_{y}} \ar[d]  & \tilde{V}_{f(y)} \ar[d] \\
\tilde{U}_y/\Gamma_y  \ar[r]^{\tilde{f}_{y}/ \Theta_{f,y}} \ar[d]    &\tilde{V}_{f(y)}/\Theta_{f,y}(\Gamma_y) \ar[d]\\
\ar[d]  &  \tilde{V}_{f(y)}/\Gamma_{f(y)} \ar[d]\\
{U}_y \ar[r]^f  & {V}_{f(y)}
 \\}
\end{equation*}

(*4) (Equivalence) Two complete orbifold maps $(f, \{\tilde{ f_x}\}, \{\Theta_{f,x}\})$ and
$(g, \{\tilde{ g_x}\}, \{\Theta_{g,x}\})$ are considered equivalent if for each $x \in X_{\mathcal O}\;\;\;   \tilde{f}_x = \tilde{g}_x $ as germs at $x$ and $\Theta_{f,x} =\Theta_{g,x}$
That is, there exists an orbifold chart $(\tilde{U}_x, \Gamma_x)$ at $x$ such 
that $\tilde{f}_x\vert _{\tilde{U}_x}  = \tilde{g}_x\vert _{\tilde{U}_x}$
and $\Theta_{f,x} =\Theta_{g,x}.$ Note that this implies that
$f = g.$

\end{definition}

\noindent
The set of smooth complete orbifold maps from $\mathcal{O}$ to $\mathcal{P}$ will be denoted by   $C^{\infty}_{Orb}(\mathcal{O}, \mathcal{P})$

\noindent

In a very similar way to submanifolds one can define suborbifolds. The first definition of a suborbifold was given by W. Thurston in \cite{Th}, but it seems to be too restrictive. We recall and use the definition formulated   in \cite{BB}

\begin{definition} An (embedded) suborbifold P of an orbifold O consists of the following:

i) A subspace $X_P \subset X_O$ equipped with the subspace topology,

ii)  For each $x \in X_P$ and neighbourhood W of x in $X_O$ there is an orbifold chart $(\tilde{U}_x, \Gamma_x, \rho_x, \phi_x)$ about x in O with 
$U_x \subset W,$ a subgroup $\Lambda_x \subset \Gamma_x$ of the
isotropy group of x in O and a $\rho_x(\Lambda_x) $-invariant submanifold $\tilde{V}_x  \subset \tilde{U}_x \cong R^n,$  so that 
$( \tilde{V}_x, \Lambda_x / \Omega_x, \rho_x\vert_{\Lambda_x}, \Psi_x ) $ is an orbifold chart for $P$ where  $\Omega_x  = \{ \gamma \ \Lambda_x \colon \rho_x(\gamma ) \vert _{\tilde{V}_x} = id \}. $ (In particular, the intrinsic isotropy subgroup at  
x), 

iii) $V_x =  \psi_x( \tilde{V}_x/\rho_x(\Lambda_x)) = U_x \cap X_P $  is an orbifold chart domain  for $x\in P.$

\end{definition}

 Thurston's definition is a bit more restrictive:

\begin{definition} A $ Y \subset X$  is called a full suborbifold of X if Y is a suborbifold with $\Lambda_x = \Gamma_x$
for all  points x of  Y.

\end{definition}

For the differences in the definitions and the reasons for them see \cite{BB}

\section{Fibre bundles over orbifolds}

Let $F$ be a smooth manifold. 

\begin{definition}
An orbifold $E$ is called an orbifold  fibre bundle over the orbifold $X$ with standard fibre $F$ if 

i) there exists a smooth orbifold map $p\colon E  \rightarrow X,$

ii) there exists an orbifold atlas $\mathcal U$  of $X,$ i.e., for any $x\in X$ there exists an orbifold chart  $(U_i, \Gamma_i, \psi_i)$ of $\; \mathcal U$  such that  $x\in U_i$  and 
 $\tilde{U_i}$ an open subset of $R^n$, $\Gamma_i$ is a finite group  of diffeomorphisms of  $\tilde{U_i},$ and $\tilde{U_i}/\Gamma_i$ is homeomorphic to $U_i$, $\psi_i$ being the homeomorphism, 

iii) let $V_i= p^{-1}(U_i)$ and   $\tilde{V_i} = \tilde{U_i} \times F$, then there exist a group $\Lambda_i$ of fibre preserving diffeomorphisms of $\tilde{V_i}$ and a  homeomorphism $\phi_i \colon \tilde{V_i} /\Lambda_i  \rightarrow V_i$  such that  $\{ ( \tilde{V}_i, \Lambda_i) \}$ form an  atlas of the orbifold $E$, 

iv) and the following diagram is commutative

\begin{equation*}
\xymatrixcolsep{7pc}\xymatrix{ 
\tilde{U_i} \times F \ar[r]^{ \tilde{p}=p\times id } \ar[d]  &  \tilde{U}_i \ar[d] \\
\tilde{V_i} /\Lambda_i  \ar[d]^{\phi_i}    & \tilde{U_i}/\Gamma_i  \ar[d]^{\psi_i}\\
{V}_i \ar[r]^p  & {U}_{i}
 \\}
\end{equation*}

\end{definition}

\medskip

\noindent
{\bf Remark} We can assume that $\tilde{p}$ is $(\Lambda_i,  \Gamma_i)$-equivariant. Obviously, $p$ is a $C^{\infty}$-complete orbifold map. In fact, we will consider only fibre budnles for which the groups $\Lambda_i$ and $\Gamma_i$ are isomorphic. 

\medskip

\noindent
{\bf Examples}

\medskip

\noindent
a)  The {\bf  tangent bundle}  $TX$ of an orbifold $X$,  $F=R^n$

\medskip

 We construct the tangent bundle as follows. Take any orbifold atlas $\;{\mathcal U} =\{ (\tilde{U}_i, \Gamma_i,  \phi_i)\}$. Then take $\tilde{V}_i = T\tilde{U}_i = \tilde{U}_i\times R^n,$ as the  group $\Sigma_i$ local transformations take 
$\Sigma_i = \{ d\gamma \colon \gamma \in \Gamma_i \}$ and as $\psi_i$ the quotient map $\tilde{V}_i \rightarrow \tilde{V}_i/\Sigma_i$. The condition of local compatibility of these charts is obviously satisfied with the required  embeddings provided by the differentials of the embeddings of the compatibility condition 
of the atlas ${\mathcal U}$. The topological space $TX$ is defined as the quotient of the disjoint union of $T\tilde{U}_i$ by the equivalence relation induced by the action of 
the groups $\Sigma_i$ and the embeddings from the definition of an orbifold atlas. Equivalent atlases of the orbifold $X$ define equivalent atlases of $TX$.

A smooth complete orbifold mapping $f \colon X \rightarrow Y$ defines the mapping $df \colon TX \rightarrow TY$ of the orbifold tangent bundles which is called the differential of $f.$ If  $(\tilde{U}, \Gamma,  \phi)$ is a chart of $X$ and $ (\tilde{W}, \Delta,  \rho)$ is a chart of $Y$ such that $f(U) \subset W$
and $\tilde{f} \colon \tilde{U} \rightarrow \tilde{W}$ covers $f$ then the differential $d\tilde{f}$ satisfies the condition of a complete orbifold mapping for the just defined atlases of the orbifold tangent bundles $TX$ and $TY$. It is a simple exercise to verify that $df$ is a smooth complete map of the orbifolds $TX$ and $TY$.  Thus we have defined  a functor from the category of (effective) orbifolds and their complete smooth mappings into the category of orbifold vector bundles and their complete smooth vector bundle mappings.

\medskip

\noindent
b) The {\bf linear frame} bundle of an orbifold, $F = GL(n)$

The same construction as in the point a) using the frame bundle of open subsets of $R^n$ instead of the tangent bundle  defines the orbifold principal frame bundle of an orbifold. As the local groups are finite, and the fact that one can make them local isometries, the total space 
of the orbifold principal bundle is in fact a manifold.  Any smooth complete mapping of orbifolds defines a smooth complete mapping of their frame bundles.   
The above correspondence  is a functor from the category of (effective) orbifolds and their complete smooth mappings into the category of principal frame  bundles and their  smooth  bundle mappings. 

\medskip
\noindent
c) {\bf Higher order frame and tangent}  bundles

\medskip

The same procedure can be applied to the functor of higher order frame and tangent  bundles. In this way we define the functor from the category of smooth orbifolds and their smooth complete maps to the category of orbifold fibre bundles and their smooth complete fibre maps.

\medskip

Similar constructions work well in the cases of fibre bundles listed below:

\medskip

    d) associated bundles to higher order frame bundles

     e)  natural bundles

\medskip

It is not difficult to verify that one can define the dual vector bundle $E^*$ of any orbifold vector bundle $E,$ and that $E^*$ is an orbifold vector bundle. Moreover, any   tensor product of  orbifold vector bundles over a given orbifold is itselt an orbifold vector bundle. Therefore we can define the  orbifold tensor algebra

$$\bigotimes TX = \bigoplus_{p,q} \bigotimes _p^q TX = \bigoplus_{p,q} (\otimes^pTX )\otimes (\otimes^q TX^*)$$

\medskip

\noindent
Likewise for any finite number of orbifold vector bundles over a given orbifold $X$ we can define their skewsymmetric product which itself is an orbifold vector bundle over $X.$

\bigskip

Another classical construction for fibre bundles over smooth manifolds works well within the framework of orbifolds, the pullback. 
Let $f \colon X \rightarrow Y$ be a  smooth complete orbifold mapping between two orbifolds $X$ and $Y$. Let $p\colon E \rightarrow Y$ be a smooth orbifold bundle, then the pullback bundle

$$f^{-1}E = \{ (x,w) \in X \times E \colon f(x) =p(w) \}$$

\noindent
is a well-defined orbifold bundle over $X$ with the same standard fibre as the bundle $E.$

\bigskip
\noindent
{\bf Remark} The theory of natural bundles is well presented in \cite{Michor}, see also \cite{terng}.

\subsection{Sections of fibre bundles over an orbifold}

We shall consider sections of orbifold fibre bundles which are smooth complete orbifold mappings.  

\medskip

If $p\colon E \rightarrow X$ is a smooth orbifold bundle, $Sect_{comp}(X,E)$ denotes the space of  smooth complete sections of $E$, i.e., the set of all smooth complete orbifold mapping $s \colon X \rightarrow E$ such that $ps =id_{X}$

Let $\mathcal U$ be an orbifold atlas of the orbifold $X.$ The existence of a section $s$ of the fibre bundle $s \colon X    \rightarrow E$ is equivalent to the existence on each $U_i$ of a  section $s_i$ of $\tilde{U}_i \times F$ which is $(\Gamma_i, \Lambda_i)$-equivariant and satisfy the compatibility condition.

\medskip

\noindent{\bf Riemannian metrics}

A Riemannian metric on an orbifold $X$ is given by a family of Riemannian metrics $g_i$ on $\tilde{U}_i$ which are $\Gamma_i$ invariant, i.e. elements of  $\Gamma_i$ are isometries of  the Riemannian metrics $g_i.$

\begin{proposition}
On any orbifold $X$  there exists a Riemannian metric $g_X. $
\end{proposition}

\medskip
As the metrics $g_i$ on $\tilde{U}_i$ are „compatible”, so are the associated Levi-Civita connections $\nabla^i$. The induced object on the orbifold $X$ we call the Levi-Civita connection of the Riemannian metric $g$. The operator $\nabla$ can be characterized as a two-linear mapping

$$ {\mathcal X}_{comp}(X) \times {\mathcal X}_{comp}(X)  \rightarrow {\mathcal X}_{comp}(X)$$

\noindent
satisfying the standard condtions for connections. 
The operator $\nabla$ can be extended to the orbifold  tensor algebra  $\bigotimes TX$ in the same way as in the case of manifolds. For a given complete vector field $Z$ on the orbifold $X$ and a complete $(p,q)$ tensor field $T, \;\;$ $\nabla_ZT$ is also a   complete $(p,q)$ tensor field, or $\nabla T$ a complete $ (p,q+1) $ tensor field.

A complete section of the orbifold tensor bundle $\otimes^p_qTX$ is called a $(p,q)$-orbi-tensor field on the orbifold $X.$   
Let $T_0 \in (\otimes^pR^n) \otimes (\otimes^qR^{n*}) $ be  a $(p,q)$-tensor of $R^n$. A tensor field $T$ is called {0-deformable\/} if for any point $x\in X$ there exists a frame $p$ at that point such that $p^*T_x = T_0$ for some tensor $T_0.$

\subsection{Geometrical structures on orbifolds}

In the case of manifolds a very general definition defines a geometrical structure as a submanifold or a subbundle of a natural bundle over the chosen manifold. The same approach can be extended to the case of orbifolds. 

A suborbifold $H$ of  the total space of an orbifold natural bundle $K$  which is at the same time an orbifold natural bundle is called a suborbifold fibre bundle. 

\begin{definition}
A geometrical structure on an orbifold X is an orbifold subbundle of a natural orbifold bundle H over X.

\end{definition}

In particular, one can develop the theory of classical orbifold G-structures for any classical Lie group $G \subset GL(n)$. For a fixed Lie group $G$  orbifold 
G-structures on an orbifold of dimension $n$ is an orbifold subbundle $B(X,G)$ of the linear frame bundle $L(X)$. As in the case of manifolds any 0-deformable orbi-tensor field $T$ determines an orbifold $G_T$-structure.  In this case the corresponding $G$-structure is defined as, for any point $x \in X,$

$$B(T)_x = \{ p \in L(X)_x\colon p^*T_x =T_0 \}$$

\noindent
and the structure group is the Lie subgroup of $GL(n)$:

$$G(T_0) = \{ A \in GL(n) \colon A^* T_0 = T_0 \}.$$

\medskip

\noindent
{\bf Remark} Many geometrical objects or procedures associated to $G$-structures over manifolds have their orbifold counterparts. It is not difficult to see that for an orbifold G-structure we can define the structure tensor in the same way as for manifolds, cf. \cite{Fu,St}. The structure tensor $c$ is a smooth function with values in

$$Hom(R^n \wedge R^n, R^n) / \partial Hom (R^n, Lie(G)).$$

\noindent
The vanishing of the structure tensor of an orbifold G-structure is equivalent to the existence of a torsionless orbifold connection in this orbifold G-structure.
Likewise, one can prove that the 1st prolongation of an orbifold G-structure is an orbifold bundle. So the theory of prolongations of G-structures works well for  orbifolds, cf. \cite{Fu,St}.

\medskip

Like in the classical case of manifolds we can define many geometrical structures in the dual way, via orbi-tensor fields or via orbifold reductions of the orbifold frame bundle. 

\medskip

\noindent
{\bf Examples} 

\medskip

a) {\bf Riemannian structure (Riemannian metric)} In the previous subsection we have defined an orbifold  Riemannian metric as a section of the orbifold bundle $\wedge^2T^*X.$ It is equivalent to a choice of an orbifold $O(n)$ reduction of the orbifold linear frame bundle $LX.$

\medskip

b) {\bf Symplectic structure } A symplectic structure on the orbifold $X$ can be introduced by choosing a symplectic form on $X,$ or equivalently an orbifold $Sp(2m)$-reduction of the orbifold linear frame bundle $L(X),$ cf. \cite{WZ,BC}.

\medskip

c) {\bf K\"ahler structure} A $U(n)$-reduction of the orbifold  linear frame bundle $LX$ of the orbifold $X$ would give us an almost K\"ahler structure on $X.$ The vanishing of its structure tensor is equivalent to the existence of a $U(n)$ torsionless connection, thus to the integrability of the almost-complex structure.

\medskip

d) In a similar way we can define on hyperK\"ahler or  quaternionic structures on orbifolds, cf. \cite{BG}. 

\medskip

{\bf Remark}
Having said that we realize that with no difficulty we can define on orbifolds contact structures, K-contact and Sasakian structures as well as 3-Sasakian or K-structures. The only difference is that all the objects considered have to be orbifold tensor fields.

\medskip
\subsection{Harmonic mappings}

Any smooth complete orbifold mapping $f \colon X \rightarrow Y$ defines a smooth complete mapping $df \colon TX \rightarrow TY$ called the differential of $f$.  The differential $df$ is a smooth complete section of the bundle

$$TX\otimes f^{-1}TY$$

\medskip

If the orbifolds $X$ and $Y$ are Riemannian, then the associated Levi-Civita connections $\nabla^X$ and $\nabla^Y$, respectively, define a connection $D$ in the vector bundle 
$TX\otimes f^{-1}TY$. Therefore it makes sense to consider

$$Ddf$$

\noindent
and its trace 

$$\tau (f)$$

\noindent
$\tau (f)$ is a complete section of the bundle $ f^{-1}TY \rightarrow X$. We call it the tension field of the complete mapping $f$  between the Riemannian orbifolds $(X,g_X)$ and $(Y,g_Y)$. 

\begin{definition}
A complete orbifold mapping $f\colon X \rightarrow Y$ between two Riemannian orbifolds $(X,g_X)$ and $(Y,g_Y)$ is harmonic if its tension field vanishes.
\end{definition}

\medskip

\noindent
{\bf Remark} A different approach to the study of harmonic maps between orbifolds was proposed by Y.-J. Chiang in \cite{C}. The same approach can be applied to the study of bi-harmonic, f-harmonic, etc. mappings  as well as harmonic morphisms between orbifolds using the theory developped by Y.-J. Chiang and the author, cf. \cite{CW1,CW2,CW4,CW5,CW6}.

\section{Orbifolds and foliations}

At the begining of the seventies foliations with all leaves compact were an object of very interesting studies. In the case of compact manifolds it was proved that the local boundedness of 
the volume of leaves is equivalent to the leaf space being Hausdorff, or to the fact that the holonomy group of any leaf  is finite. That in turn means that the leaf space of such a foliation is 
an orbifold, cf. \cite{Eps,Mi,EMS}. In this case the foliation is Riemannian and taut, cf.  \cite{Molino}. A few years later A. Haefliger and J. Girbau, cf.  \cite{HG}, remarked that any orbifold can be realized as the leaf space of a Riemannian
foliation with compact leaves. It is a simple remark, having noticed that an orbifold admits a Riemannian metric, they take the associated orthonormal frame bundle, which in this case is a smooth manifold. If we take a compact orbifold, 
its orthonormal frame bundle is a compact manifold foliated by a foliation with all leaves compact, the fibres of this bundle. Its leaf space can be identified with    the initial orbifold. Therefore Riemannian foliations with all leaves compact, or more general 
Riemannian foliations with leaves of finite holonomy, can be considered as natural desingularization objects of orbifolds. For a given orbifold $X$ we will denote this foliated manifold by $(O_X, {\mathcal F}_X).$

In the study of the geometry of a single orbifold this duality between the orbifold and its orthogonal frame bundle foliated by the compact fibres, between orbifold objects and foliated objects, works well, cf. e.g. \cite{BC}.

The theory of foliated G-structures was developed mainly by Pierre Molino, cf. \cite{Molino} and his numerous papers. The most general theory was proposed by the author in \cite{Wo0}. It is based on the theory of natural bundles.


Let $(M,g,{\mathcal F})$ be a compact Riemannian foliated manifold with compact leaves. In fact, we can assume less: it is sufficient to have a complete bundle-like  Riemannian metric and leaves of finite holonomy. Then according to  the  Reeb stability theorem, cf. \cite{Re,CN},  each leaf $L$ admits a foliated (saturated by leaves) tubular neighbourhood $U.$ Any fibre of such a tubular neighbourhood is a transverse submanifold intersecting all leaves contained in the neighbourhood. The holonomy of the zero section, a leaf of $\mathcal F$, defines a finite subgroup $\Gamma$ of diffeomorphisms of the fibre $D.$ The space $U/{\mathcal F}$ of leaves  contained in the tubular neighbourhood $U$ is identified with the the orbit space $D/\Gamma .$ Let $U$ and $V$ be two such tubular neighbourhoods of leaves $L_1$ and $L_2,$ respectively.  If the intersection is not empty, it is an open saturated subset. For any leaf $L$ in $U\cap V$ we can find an open saturated  tubular neighbourhood $W \subset U \cap V.$ The choice of a fibre of $W$ defines two embeddings of $W/{\mathcal F}$ into $U/{\mathcal F}$  and $V/{\mathcal F},$  respectively. Thus the orbifold charts of the leaf space $M/{\mathcal F}$ defined using these tubular neighbourhoods are compatible, hence we have an orbifold atlas on $M/{\mathcal F}.$  Different choices of tubular neighbourhoods and of different fibres in them produce different but equivalent orbifold atlases. Thus  we have defined an orbifold structure on $M/{\mathcal F}.$

In the language of foliations, the  orbifold associated to a foliations with compact leaves and finite holonomy contains the following data: the underlying topological space is the  space of leaves of the foliation, the orbifold atlas corresponds to the choice of a complete transverse manifold, 
and the local finite groups are the elements of the holonomy pseudogroup induced on this transverse manifold. The embeddings from the definition of the orbifold atlas are also induced by elements of the holonomy pseudogroup, those which map open subsets of one connected component of the transverse manifold into another. The choices described in the previous paragraph lead to different but equivalent pseudogroups, cf. \cite{Haef}.

In \cite{Wo0},  the author demonstrated that foliated geometrical structures are in one-to-one correspondence with holonomy invariant geometrical structures on the transverse manifold.  In the case of foliated manifolds with foliations with compact leaves and finite holonomy, this result can be reformulated saying  that foliated geometrical structures are in one-to-one correspondence with orbifold geometrical structures on the leaf space equiped with the induced orbifold structure. In particular, any bundle-like Riemannian metric of the foliated manifold $(M,{\mathcal F})$ induces 
an orbifold Riemannian metric on the leaf space  $M/{\mathcal F}$ and vice versa.

\medskip

Let $f\colon (M_1,{\mathcal F}_1) \rightarrow (M_2,{\mathcal F}_2) $ be a foliated smooth mapping between two foliated manifolds. We can choose two open coverings ${\mathcal U}_1$ and ${\mathcal U}_2$ by tubular neighbourhoods of leaves of 
$(M_1,{\mathcal F}_1)$  and $(M_2,{\mathcal F}_2) $, respectively, such that for any $U \in {\mathcal U}_1$ there exists $V \in {\mathcal U}_2$  such that $f(U) \subset V$. This choice ensures that the foliated mapping $f$ induces a holonomy invariant mapping $\bar{f}$  of the associated transverse manifolds,
cf. \cite{KW1}. Therefore the induced  mapping $\tilde{f} \colon M_1/{\mathcal F}_1 \rightarrow M_2/{\mathcal F}_2$ between the leaf spaces, considered as orbifolds, is a complete orbifold mapping. The remark can be understood as an orbifold version of the considerations of Section 3.2 of \cite{KW1} about the mapping induced on the transverse manifolds by a foliated mapping.

\medskip

Let us recall the precise definition of the mapping $\bar{f}$ presented in \cite{KW1}.

\medskip

Let $f\colon (M_1,{\mathcal F}_1) \rightarrow (M_2,{\mathcal F}_2) $ be a foliated smooth mapping between two foliated manifolds. 
 We suppose that ${\mathcal U}  =  \{(U_i, \phi_i , g_{i j })\} _I$  is a cocycle defining the foliation ${\mathcal F}_1,$  and denote by 
 ${\mathcal V}  =  \{(V_{\alpha}, \psi_{\alpha} , h_{\alpha\beta })\} _A$  a cocycle defining the foliation
${\mathcal F}_2$ such that, for any  $i \in I ,$ there exists $\alpha (i ) \in A$  such that $f (U_i ) \subset  V_{\alpha (i )}.$ Let $\bar{U}_i = \phi_i (U_i )$ and
$\bar{V}_{\alpha} = \psi_{\alpha} (V_{\alpha} ).$ Then the manifold $N_1 = \coprod \bar{U}_i $ is a transverse manifold of the foliation ${\mathcal F}_1$, and
$N_2 = \coprod \bar{V}_{\alpha} $ is a transverse manifold of the foliation ${\mathcal F}_2 .$ The transformations $g_{i j}$ generate a
pseudogroup ${\mathcal H}_1, $ which is called the holonomy pseudogroup of ${\mathcal F}_1$ associated with the cocycle ${\mathcal U},$
and the transformations $h_{\alpha\beta}$ generate a pseudogroup ${\mathcal H}_2, $ which is called the holonomy pseudogroup
of ${\mathcal F}_1$ associated with the cocycle ${\mathcal V}.$ On the level of transverse manifolds, the map $f $ induces a smooth map $\bar{f},$  as, for any $i \in I ,$  the following diagram is commutative:

\begin{equation*}
\xymatrixcolsep{7pc}\xymatrix{ 
U_i  \ar[r]^{f\vert U_i} \ar[d]^{\tilde{\phi}_{i}  } &  \tilde{U}_x \ar[d]^{\tilde{\psi}_{\alpha (i)}} \\
\bar{U}_i \ar[r]^{\bar{f}_{\alpha (i)i}}  & \bar{V}_{\alpha (i)}
 \\}
\end{equation*}

\noindent
The map $\bar{f} \colon N_1 \rightarrow N_2$ is defined as follows:

$$ \bar{f}\vert \bar{U}_i = \bar{f}_{\alpha (i)i}.$$

If we take two open subsets $U_i$ and $U_j$ with nonempty intersection, then $f(U_i\cap U_j) \subset V_{\alpha (i)} \cap V_{\alpha (j)}.$ This intersection covers the open subsets $\bar{U}_{ji} \subset \bar{U}_i$ and $\bar{U}_{ij} \subset \bar{U}_j.$ Likewise $V_{\alpha (i)} \cap V_{\alpha (j)}$ covers 
$\bar{V}_{\alpha (j) \alpha (i) } \subset \bar{V}_{\alpha (i)}$ and $\bar{V}_{\alpha (i) \alpha (j) } \subset \bar{V}_{\alpha (j)}.$ 
Moreover, the map 
$g_{ji} \colon \bar{U}_{ji} \rightarrow \bar{U}_{ij}$ 
is a diffeomorphism as well as the map 
$h_{\alpha (j) \alpha (i)}  \colon \bar{V}_{\alpha (j) \alpha (i) }\rightarrow \bar{V}_{\alpha (i) \alpha (j) }. $ Then

$$h_{\alpha (j) \alpha (i) }\bar{f}_{\alpha (i) i }\vert \bar{U}_{ji} = \bar{ f}_{\alpha (j) j } g_{ji} \vert \bar{U}_{ji}.$$

\medskip

Let $X$ and $Y$ be two orbifolds, and $O(X)$ and $O(Y)$ their respective foliated desingularizations, i.e., foliated manifolds whose leaf space is $X$ and $Y$, respectively. It would have been very nice if any smooth complete orbifold mapping $f \colon X \rightarrow Y$  could be lifted to a smooth foliated mapping $O(f) \colon O(X) \rightarrow O(Y).$ 
 In general, this is not the case as the differential $df$ does not send orthonormal frames into orthonormal frames. It is the case only if $f$ is an isometry of $(X,g_X)$ into $(Y,g_Y).$ To rectify this we can consider linear frame bundles $L(X)$ and $L(Y)$ of the orbifolds $X$ and $Y,$ respectively. 
Any smooth complete orbifold  local diffeomorphism $f \colon X \rightarrow Y$ lifts to a smooth mapping $L(f) \colon L(X) \rightarrow L(Y)$ which is a foliated mapping for the natural foliations ${\mathcal F}_X$ and ${\mathcal F}_Y$ by fibres of these two orbifold fibre bundles$L(X)$ and $L(Y),$ respectively. These foliations are Riemannian for the natural liftings of the Riemannian metrics of orbifolds. The holonomy groups of their  leaves are finite.

\section{Transversely harmonic maps}

Let $\mathcal F$ be a foliation on a Riemannian $n$-manifold $(M,\, g).$ Then $\mathcal F$ is defined by a cocycle    $\mathcal U =\{ U_i , f_i , g_{ij} \}_{i,j\in I}$ modeled on a $q$-manifold $N_0$ such that 

(1) $\{ U_i \}_{i\in I}$ is an open covering of M, 

(2) $f_i : U_i \to N_0$ are submersions with connected fibres, 

(3) $g_{ij} : N_0\to N_0 $ are local diffeomorphisms of $N_0$ with $f_i = g_{ij} f_j$ on $U_i \cap U_j .$ 

\noindent
The connected components of the trace of any leaf of $\mathcal F$ on $U_i$ consist of the fibres of $f_i .$ The open subsets
$N_i = f_i ( U_i )\subset N_0$ form a $q$-manifold $N=\amalg N_i$,
which can be considered as a transverse manifold of the
foliation $\mathcal F .$ The pseudogroup $\mathcal H_N$ of local
diffeomorphisms of $N$ generated by $g_{ij}$ is called the holonomy
pseudogroup of the foliated manifold $(M, \mathcal F )$ defined by the
cocycle $\mathcal U .$ If the foliation $\mathcal F$ is Riemannian for the
Riemannian metric $g,$ then it induces a Riemannian metric $\bar{g}$
on $N$ such that the submersions $f_i$ are Riemannian submersions
and the elements of the holonomy group are local isometries.

   Let  $\phi : U\to R^p \times R^q ,\,\phi =(\phi^1 ,\phi^2)=(x_1, ..., x_p ,y_1,...,y_q )$ be an adapted chart on a
   foliated manifold $(M, \mathcal F)$. Then on $U$   the  vector fields
   $\frac{\partial}{\partial x_1},...\frac{\partial}{\partial x_p}$ span
   the bundle $T\mathcal F$ tangent to the leaves of the foliation $\mathcal F ,$ the equivalence
   classes denoted by $\bar{\frac{\partial}{\partial y_1}}, ...\bar
   {\frac{\partial}{\partial y_q}}$ of  $\frac{\partial}{\partial y_1}, ...
   \frac{\partial}{\partial y_q}$ span the normal bundle $N(M, \mathcal F )
   =TM / T\mathcal F $ which is isomorphic to the subbundle $T\mathcal F^\perp$. 
 
Suppose that $(M, \, \mathcal F, \, g)$ is a Riemannian foliation. The sheaf $\Gamma_b(T{\mathcal F}^\perp )$ of foliated sections of the vector bundle ${T\mathcal F}^\perp \to M$ may be described as follows: If $U$ is an open subset of $M$, then $X\in \Gamma_b (U,T{\mathcal F}^\perp)$ if and only if for each local Riemannian submersion
$\phi :U\to \bar U$ defining $\mathcal F$, the restriction of $X$ to $U$ is projectable via
the map $\phi$ on a vector field $\bar X$ on $\bar U .$

The vector bundle $N(M, {\mathcal F})$ is a foliated vector budnle over $(M,{\mathcal F}),$ i.e., it admits a foliation ${\mathcal F}_N$ of the same dimension as ${\mathcal F}$ whose leaves are covering spaces of leaves of ${\mathcal F}.$ The vector bundle ${\mathcal F}^\perp$ is also foliated as it is isomorphic to the normal bundle. As ${\mathcal F}^\perp$ is foliated, so is its dual bundle ${\mathcal F}^{\perp *}$ and any tensor product of these vector bundles as well as any pull-back of a foliated vector bundle by a foliated mapping. 

\medskip

\begin{definition} A $basic\, partial\,
connection$ $(M,\mathcal F ,g)$ is a sheaf operator 

\noindent $D:\Gamma_b (U,T{\mathcal F}^\perp )\times  \Gamma_b (U,T{\mathcal F}^\perp )\to \Gamma_b (U,T{\mathcal F}^\perp )$ such that 
 
1. $D_{fX+hY} Z\,=\,fD_X Y+h D_Y Z ,$ 

2. $D_X$ is R-linear, 

3. $D_X fY\,=\,X(f)Y+fD_X Y $ (the transversal Leibniz rule),

\noindent for any $X, Y, Z \in \Gamma_b (U,T\mathcal F^\perp )$ and any basic functions $f,h \in C_b^\infty (U)$, where $U$ is any open subset of $M$.

\end{definition}

Let $\nabla$ be the Levi-Civita connection of g, then for any open subset U of M and $X, Y\in\Gamma_b (U,T{\mathcal F}^\perp )$ we define D as

$$ D_X Y=(\nabla_X Y )^\perp $$

\noindent
where ${(\nabla_X Y)}^\perp$ is the horizontal component of ${(\nabla_X Y)}.$ It is easy to check that $D$ is a basic partial connection on $(M,\,\mathcal F,\, g) .$ Let $\phi :U\to \bar U$ be a Riemannian submersion defining the foliation $\mathcal F$ on an open set $U.$ Let us assume that
 $X,Y\in \Gamma_b (U,T{\mathcal F}^\perp )$, and $\bar X , \bar Y$ be the push forward vector fields via the map $\phi .$ Then 

$$ d\phi (D_X Y )=\nabla_{\bar X}^{\bar g} \bar Y $$

\noindent
 where $\nabla^{\bar g}$ is the Levi-Civita connection of the metric $\bar g $ on the transverse manifold.

\medskip

Let $f:(M_1 ,{\mathcal F}_1, g)\to (M_2 ,{\mathcal F}_2 , h)$ be a foliated map between two foliated Riemannian manifolds.  The covariant derivative $D(\Pi_2dfi_1)$ is a global foliated section of the bundle

$$(T{\mathcal F}^{\perp}_1)^* \otimes (T{\mathcal F}^{\perp}_1)^* \otimes f^{-1}T{\mathcal F}^{\perp}_2 \rightarrow M_1$$

\noindent
where $i_1$ is the inclusion of 
$T{\mathcal F}^{\perp *}_1 $ into $TM_1$ and $\Pi_2$ is the orthogonal projection of $TM_2$ onto $T{\mathcal F}^{\perp}_2.$  Its trace 

$$\tau _b (f)$$

\noindent
is a foliated section of the bundle $ f^{-1}T{\mathcal F}_2^{\perp} \rightarrow M_1$. We call it the tension field of the foliated mapping $f$. 

\begin{definition}

  A foliated map $f:(M_1 ,{\mathcal F}_1, g)\to (M_2 ,{\mathcal F}_2 , h)$ between two foliated Riemannian manifolds is called transversally harmonic if
$\tau_b(f) =0$.

\end{definition}

\bigskip

As we have shown any foliated map  $f:(M_1 ,{\mathcal F}_1, g)\to (M_2 ,{\mathcal F}_2 , h)$ induces a map $\bar{f}$ of the corresponding transverse manifolds. Since the tension tensor of a foliated map is itself a foliated section, it induces a holonomy invariant object on the transverse manifold. In fact, if $p_i$ is a local submersion defining the foliation ${\mathcal F}_i$ for $i=1,2,$ then we have the following relation between the tension tensors of the mappings $f$ and $\bar{f}.$

$$dp_2\tau_b(f)_x = \tau (\bar{f})_{p_1(x)} $$

\medskip

 Let $f:(M_1 ,{\mathcal F}_1, g)\to (M_2 ,{\mathcal F}_2 , h)$ be a foliated mapping. The restriction of $f$ to any leaf $L$ of ${\mathcal F }_1$ is a smooth mapping of the manifold $L$ into a leaf of the foliation ${\mathcal F}_2$. Both submanifolds can be provided with the restrictions of the Riemannian metrics $g$ and $h,$ respectively. If for each leaf $L$ of the foliations ${\mathcal F}_1$ the mappings just described are harmonic, we say that the mapping $f$ is {\it leaf-wise harmonic.}

\medskip

Taking into account our previous considerations we can formulate the following theorems, whose proofs can be found in \cite{KW1,KW2}:

\begin{theorem} Let $f:(M_1 ,{\mathcal F}_1, g)\to (M_2 ,{\mathcal F}_2 , h)$ be a smooth foliated mapping between regular
foliated Riemannian manifolds. Let $\mathcal U$ and $ \mathcal V$ be two cocycles defining the foliations ${\mathcal F}_1$ and ${\mathcal F}_2,$  repectively. 
Then the  map $ f$  is transversally harmonic if and only if the induced map $\bar{f} \colon N_{\mathcal U} \rightarrow N_{\mathcal V}$ between the associated transverse manifolds is harmonic.
\end{theorem}


\begin{theorem} Let $f:(M_1 ,{\mathcal F}_1, g)\to (M_2 ,{\mathcal F}_2 , h)$  be a foliated harmonic map between two manifolds with Riemannian foliations. Moreover, assume that all leaves of the foliation ${\mathcal F}_1$
are minimal, the foliation ${\mathcal F}_2$ is totally geodesic, and $ f$ is horizontal, that is, for all $d f (T{\mathcal F}_1^{\perp} ) \subset  T{\mathcal F}_2^{\perp},$ then the map $f$ is transversally harmonic.
\end{theorem}


\begin{corollary} If the foliation ${\mathcal F}_1$ is minimal, and ${\mathcal F}_2$ totally geodesic, then the map $f$ is harmonic if and only if $f$ is transversally harmonic and  leaf-wise harmonic.
\end{corollary}

\begin{theorem} Let $f:(M_1 ,{\mathcal F}_1, g)\to (M_2 ,{\mathcal F}_2 , h)$ be a smooth foliated mapping between 
foliated Riemannian manifolds. If the foliations ${\mathcal F}_1$ and ${\mathcal F}_2$ are totally geodesic and the horizontal distribution of ${\mathcal F}_2$ is integrable, then the map $ f$ is transversally harmonic if and only if  the horizontal part of the tension tensor $ \tau (f )$ is zero. In particular, the harmonicity of $f$ implies its transverse harmonicity.
\end{theorem}

The above results are  "foliated" versions of Theorem 1.3.5 and its corollary of Smith's Ph.D. thesis, cf. \cite{Sm}, pp. 17 and 18, and the results of  Xin, Theorem 6.4 and its corollary of \cite{Xb}, see also \cite{X}.

\medskip

Finally, let us  formulate the following two  theorems for maps between Riemannian orbifolds which can serve as the starting point  for the study of harmonic maps between orbifolds. The proofs are the simple consequence of our previous considerations.

\begin{theorem} Let $f:(M_1 ,{\mathcal F}_1, g)\to (M_2 ,{\mathcal F}_2 , h)$ be a foliated map between two foliated Riemannian manifolds with compact leaves. Let $\bar{f}$ be the induced map between the leaf spaces with induced  orbifold Riemannian metrics, 
$\bar{f} \colon M_1 /{\mathcal F}_1, \to  M_2 /{\mathcal F}_2 .$ The map $f$ is transversely harmonic iff the map $\bar{f}$ is harmonic.
\end{theorem}

\begin{theorem} Let $f\colon X \to Y$ be a smooth complete embedding of a  Riemannian orbifolds $(X,g)$ into another Riemannian orbifold $(Y,h).$ Then $f$ is harmonic iff the induced (foliated) mapping $L(f) \colon (L(X), {\mathcal F}_X, g_L)  \to (L(Y), {\mathcal F}_Y, h_L)$ is transversely harmonic.
\end{theorem}

\medskip

\noindent
{\bf Remark} The proposed approach the study of harmonic maps between orbifolds has some deficiences. We can not lift any smooth map to a foliated map of its foliated desingularizations. However, if we just want to find harmonic maps, then the proposed methods can offer some solutions. The foliations ${\mathcal F}_X$ and ${\mathcal F}_Y$ are taut, so harmonic foliated maps between $O_X$ and $O_Y$ are not far from being transversally harmonic, and thus from  inducing harmonic maps between the orbifolds. On the other hand, any compact Sasakian manifold fibers over a K\"ahlerian orbifold, cf. \cite{BG},  so the study of geometrical and cohomological properties of orbifolds can lead to the development of new obstructions to the existence of such structures.

\end{document}